\documentclass[twocolumn,noversion,nonote,squeezed]{cdcarticle}

\def\MODE{3}

\usepackage{cite}
\usepackage{math}						      
\usepackage[utf8]{inputenc}				
\usepackage[english]{babel}				
\usepackage[hidelinks]{hyperref}	
\usepackage{enumitem}					    
\usepackage{amsmath,amssymb,amsfonts}
\usepackage{tablefootnote}
\usepackage{textcomp}
\def\BibTeX{{\rm B\kern-.05em{\sc i\kern-.025em b}\kern-.08em
    T\kern-.1667em\lower.7ex\hbox{E}\kern-.125emX}}

\newcommand{\1}{\mathbf{1}}

\newcommand{\df}{\nabla\! f}

\usepackage{algorithm}
\floatname{algorithm}{Canonical Form}

\usepackage{algpseudocode}
\usepackage{xpatch}
\def\mathllap{\mathpalette\mathllapinternal}
\def\mathllapinternal#1#2{\llap{$\mathsurround=0pt#1{#2}$}}
\algrenewcommand\algorithmicindent{0.8em}  
\makeatletter
\xpatchcmd{\algorithmic}{\itemsep\z@}{\itemsep=0.5ex plus2pt}{}{}
\makeatother

\let\oldbibliography\thebibliography
\renewcommand{\thebibliography}[1]{\oldbibliography{#1}
\setlength{\itemsep}{1.6pt}} 

\begin{document}

\title{A Canonical Form for First-Order\\ Distributed Optimization Algorithms}

\if\MODE1
\author{Akhil Sundararajan \and Bryan Van Scoy \and Laurent Lessard}
\else
\author{Akhil Sundararajan\footnotemark[1] \and Bryan Van Scoy\footnotemark[2] \and Laurent Lessard$^{1,2}$}\fi

\note{American Control Conference, pp. 4075--4080, 2019.}
\maketitle



\if\MODE1\else

\footnotetext[1]{A.~Sundararajan and L.~Lessard are with the Department of Electrical and Computer Engineering at the University of Wisconsin--Madison, Madison, WI~53706, USA.}
\footnotetext[2]{B.~Van Scoy and L.~Lessard are with the Wisconsin Institute for Discovery, which is also at the University of Wisconsin--Madison.\\
\texttt{\{asundararaja,vanscoy,laurent.lessard\}@wisc.edu}\\[1mm]
This material is based upon work supported by the National Science Foundation under Grant No. 1656951.}

\fi

\begin{abstract}
We consider the distributed optimization problem in which a network of agents aims to minimize the average of local functions. To solve this problem, several algorithms have recently been proposed where agents perform various combinations of communication with neighbors, local gradient computations, and updates to local state variables. In this paper, we present a canonical form that characterizes any first-order distributed algorithm that can be implemented using a single round of communication and gradient computation per iteration, and where each agent stores up to two state variables. The canonical form features a minimal set of parameters that are both unique and expressive enough to capture any distributed algorithm in this class. The generic nature of our canonical form enables the systematic analysis and design of distributed optimization algorithms.
\end{abstract}

%

\section{Introduction}\label{sec:intro}

We consider the distributed optimization problem in which a network of $n$ agents are connected through an undirected graph. The global objective is to minimize the average of $n$ functions $f_i:\R^d\rightarrow \R$ that are local to each agent $i=1,\dots, n$.  While each agent has access to its local variable $x_i$, it must exchange information with its neighbors to arrive at the global minimizer
\begin{equation}\label{eq:distrop-problem}
x^\star \defeq \argmin_{x\in \R^d} f(x),
\quad\text{where }f(x) \defeq \frac{1}{n}\sum_{i=1}^n f_i(x).
\end{equation}
Distributed optimization is an active area of research and has attracted significant attention in recent years due to its applications in distributed machine learning, distributed estimation, and resource allocation~\cite{forero2010consensus,johansson2008distributed,predd2009collaborative,ling2010decentralized}.  To solve~\eqref{eq:distrop-problem}, many gradient-based algorithms have been proposed involving gradient computation steps in conjunction with gossip steps, i.e., local averaging operations~\cite{DGD, EXTRA, NIDS, DIGing, ExactDiffusion, AugDGM}. Such algorithms have much structural variety in how they perform gossip, evaluate the gradient, and update their local state variables. Furthermore, changing the order of these communication and computation steps results in different algorithms. As an example, consider the NIDS~\cite{NIDS} algorithm update:
\begin{align*}
x_i^0 &\in\R^d\text{ arbitrary,}\quad x_i^1 = x_i^0 - \alpha \df_i(x_i^0) \tag{NIDS}\\
x_i^{k+2} &= \sum_{j=1}^n \widetilde{w}_{ij} \bigl(2 x_j^{k+1} - x_j^k - \alpha\df_j(x_j^{k+1}) + \alpha\df_j(x_j^k)\bigr)
\end{align*}
where $\{\widetilde{w}_{ij}\}\in\R^{n\times n}$ is a gossip matrix, $\alpha\in\R$ the stepsize, and $x_i^k\in\R^d$ the state of agent~$i$ at iteration~$k$. Each update involves two previous iterates as well as a difference of gradients, similar to so-called ``gradient tracking" algorithms~\cite{DIGing}. In contrast, Exact Diffusion~\cite{ExactDiffusion} uses two state variables, which are updated with a gradient step followed by a gossip step:
\begin{equation}\label{eq:ExactDiffusion}\tag{Exact Diffusion}
\begin{aligned}
x_{1i}^{k+1} &= x_{2i}^{k}-\alpha\,\df_i (x_{2i}^k) \\ 
x_{2i}^{k+1} &= \sum_{j=1}^n w_{ij}\,(x_{1j}^{k+1} - x_{1j}^k + x_{2j}^k).
\end{aligned}
\end{equation}
While their update equations appear different, we show in Section~\ref{sec:algorithms} that NIDS and Exact Diffusion are actually equivalent, which exemplifies the need to better understand the taxonomy of distributed algorithms.  To this end, we develop a canonical form that parameterizes the following broad class of first-order distributed algorithms.
\paragraph{Algorithm properties.}\label{alg_class} \textit{We consider the class of distributed algorithms that satisfy the following properties.}
\begin{enumerate}[itemsep=0mm,topsep=1mm,label=P\arabic*.,ref=\textup{P\arabic*}]
	\item \textit{Each agent has a local state of dimension at most~$2d$. Since $\df_i:\R^d\to\R$ by assumption, this allows each agent to store up to two past iterates.}\label{P1}
	\item \textit{At each iteration, each agent may do each of the following once in any order:}\label{P2}
	\begin{enumerate}[itemsep=0mm,topsep=0mm,label=(\alph*)]
		\item \textit{communicate any number of local variables with its immediate neighbors simultaneously,}
		\item \textit{compute its local gradient at a single point, and}
		\item \textit{update its local state.}
	\end{enumerate}
	\item \textit{The local state updates are linear, time-invariant, deterministic, and homogeneous across all agents and dimensions of the objective function.}\label{P3}
\end{enumerate}
While NIDS and Exact Diffusion belong to the class of algorithms described by \ref{P1}--\ref{P3}, not all algorithms do. Exceptions include accelerated~\cite{QuLi_accelerated,DHB}, proximal~\cite{PG-EXTRA}, stochastic~\cite{DGST}, and asynchronous~\cite{AsynDGM} variants.

\subsection{Canonical form}

Similar to the controllable canonical form for linear time-invariant systems, we want our canonical form to have existence and uniqueness properties.  That is, any algorithm satisfying \ref{P1}--\ref{P3} should be equivalent to the canonical form through appropriate selection of parameters, with a one-to-one correspondence between algorithms and parameter selections.  The canonical form should also be minimal in that it uses the fewest number of parameters possible to express all algorithms in the considered class.

In the canonical form, we use the graph Laplacian matrix $L\in\R^{n\times n}$ to model local communication. This is equivalent to using a gossip matrix $W$ with $W\defeq I-L$. We make the following assumption throughout the paper.

\begin{assumption}\label{assumption:graph} The Laplacian matrix $L\in \R^{n\times n}$ is symmetric, positive semidefinite, and satisfies $L\1=0$, where $\1$ is the all-ones vector.  The zero eigenvalue of~$L$ is distinct, implying that the underlying graph is connected.
\end{assumption}
\noindent
Our canonical form satisfies properties \ref{P1}--\ref{P3}, and is parameterized by five real scalars: $\alpha,\zeta_0,\zeta_1,\zeta_2,\zeta_3$.  
\vspace{-3mm}
\begin{algorithm}[H]
	\small\caption{}\label{alg:canonical}
	\begin{algorithmic}
		\State\textbf{Initialization:}~Let $L\in\R^{n\times n}$ be a Laplacian matrix. Each agent $i\in\{1,\ldots,n\}$ chooses $x_i^0 \in \R^d$ and $w_i^0\in\R^d$.
		\For {iteration $k=0,1,2,\ldots$}
		\For {agent $i\in\{1,\ldots,n\}$}
		\State \textbf{Local communication}
		\State $~~\hphantom{w_i^{k+1}}\mathllap{v_{1i}^k} = \sum_{j=1}^n L_{ij}\,x_j^k \hfill \text{(C.1)}$
		\State $~~\hphantom{w_i^{k+1}}\mathllap{v_{2i}^k} = \sum_{j=1}^n L_{ij}\,w_j^k \quad \text{(only required if } \zeta_2 \ne 0) \hfill \text{(C.2)}$
		\State \textbf{Local gradient computation}
		\State $~~\hphantom{w_i^{k+1}}\mathllap{y_i^k} = x_i^k - \zeta_3\,v_{1i}^k \hfill \text{(C.3)}$
		\State $~~\hphantom{w_i^{k+1}}\mathllap{u_i^k} = \df_i(y_i^k) \hfill \text{(C.4)}$    
		\State \textbf{Local state update}
		\State $~~\hphantom{w_i^{k+1}}\mathllap{x_i^{k+1}} = x_i^k + \zeta_0\,w_i^k - \alpha\,u_i^k - \zeta_1 v_{1i}^k + \zeta_2 v_{2i}^k \hfill \text{(C.5)}$
		\State $~~\hphantom{w_i^{k+1}}\mathllap{w_i^{k+1}} = w_i^k - v_{1i}^k \hfill \text{(C.6)}$
		\EndFor
		\EndFor
	\end{algorithmic}
\end{algorithm}
\vspace{-4mm}
The algorithm requires agent $i$ to store two local state variables, $x_i^k$ and $w_i^k$. Agents first communicate their states with neighbors using the Laplacian matrix; note that $w_i^k$ need not be transmitted if $\zeta_2=0$ since the result is unused in that case. Agent $i$ then evaluates its local gradient at $y_i^k$, and updates its states using the results of the communication and computation. The sequence $\{y_i^k\}$ is then agent $i$'s estimate of the optimizer~$x^\star$.

\paragraph{Technical conditions.} \textit{The parameters $\alpha,\zeta_0,\zeta_1,\zeta_2,\zeta_3$ are constant scalars that satisfy the following:}
\begin{enumerate}[itemsep=0mm,topsep=1mm,label=T\arabic*.,ref=\textup{T\arabic*}]
\item $\alpha \ne 0$.\label{T1}
\item \textit{The linear system $\alpha u = (\zeta_0 I + \zeta_2 L) w$ has a solution $w\in\R^n$ for all $u\in\R^n$ that satisfies $\1^\tp u = 0$.} \label{T2}
\item \textit{Either $\zeta_0=0$ or the algorithm is initialized such that $\sum_{i=1}^n w_{i}^0=0$. An easy way to satisfy this condition is for every agent to use the initialization $w_i^0=0$.} \label{T3}
\end{enumerate}

\subsection{Main results}

We present the two main results of the paper, which relate the distributed algorithm class to our canonical form. 

\begin{defn}[Optimal fixed point]\label{def:opt_fixed_pt}
The canonical form has an \emph{optimal fixed point} if there exists a fixed point $(v_{1i}^\star,v_{2i}^\star,y_i^\star,u_i^\star,x_i^\star,w_i^\star)$ for $i\in\{1,\dots,n\}$ such that ${y_i^\star = x^\star}$ for all $i\in\{1,\ldots,n\}$, where $x^\star$ is defined in~\eqref{eq:distrop-problem}.
\end{defn}

Our first main result states that the technical conditions are \emph{sufficient} for the canonical form to have an optimal fixed point; we prove the result in Appendix~\ref{sec:fixedpointproof}.

\begin{thm}[Sufficiency]\label{thm:fixed_point}
If technical conditions~\ref{T1}--\ref{T3} are satisfied, then the canonical form has a fixed point. Moreover, all fixed points are optimal fixed points.
\end{thm}

Note that the technical conditions are sufficient for the canonical form to have an optimal fixed point, but do \textit{not} guarantee convergence to the fixed point.
Our second main result is that the technical conditions are \textit{necessary} for convergence to an optimal fixed point. The result also characterizes the existence and uniqueness properties of our canonical form. We prove the result in Appendix~\ref{sec:realizationproof}.

\begin{thm}[Necessity]\label{thm:realization}
Consider a distributed algorithm that satisfies properties~\ref{P1}--\ref{P3} and converges to an optimal fixed point. Then there exist parameters $\alpha,\zeta_0,\zeta_1,\zeta_2,\zeta_3$ for which the algorithm is equivalent to our canonical form. Furthermore, the parameters are unique, satisfy technical conditions~\ref{T1}--\ref{T3}, and constitute a minimal parameterization of all such algorithms.
\end{thm}

To prove the main results, we develop the notion of a transfer function for distributed algorithms in Section~\ref{sec:transferfunction} and provide associated necessary conditions that hold when the algorithm converges to an optimal fixed point. We then show how to put distributed algorithms into canonical form through an example in Section~\ref{sec:algorithms}.

\section{Transfer function interpretation}\label{sec:transferfunction}

In this section, we show how to represent distributed algorithms as dynamical systems and obtain  corresponding transfer functions. We derive necessary conditions on the poles and zeros of the transfer function of any distributed algorithm that solves~\eqref{eq:distrop-problem}. These conditions are used to prove Theorem~\ref{thm:realization}, and also give rise to an impossibility result stating that no algorithm with a single state can solve the distributed optimization problem.

The general dynamical system model for a distributed optimization algorithm satisfying~\ref{P1}--\ref{P3} is given by\footnote{The iterations~\eqref{eq:alg} are similar to the polynomial linear protocol~\cite{PLP} used to study the dynamic average consensus problem.}
\begin{subequations}\label{eq:alg}
\begin{align}
\bmat{\xi_i^{k+1} \\ y_i^k} &= \biggl(\bmat{A_0 & B_0 \\ C_0 & D_0}\otimes I_d\biggr) \bmat{\xi_i^k \\ u_i^k} \nonumber\\
&\qquad + \biggl(\,\bmat{A_1 & B_1 \\ C_1 & D_1}\otimes I_d\biggr)\sum_{j=1}^n L_{ij}  \bmat{\xi_j^k \\ u_j^k} \\
u_i^k &= \df_i(y_i^k)
\end{align}
\end{subequations}
for $i\in\{1,\ldots,n\}$ where $\xi_i^k \in \R^{2d}$ is the state, $u_i^k \in \R^d$ the input, and $y_i^k \in \R^d$ the output of agent~$i$ at iteration~$k$. The model~\eqref{eq:alg} over-parameterizes all algorithms in our class of interest. To simplify notation, we set $d=1$ throughout the rest of the section, though our results hold for general $d\in\mathbb{N}$.

We can write~\eqref{eq:alg} in vectorized form by concatenating the states of the agents as $\xi^k \defeq \bmat{ (\xi_1^{k})^\tp & \cdots & (\xi_n^{k})^\tp}^\tp$ and similarly for $u^k$ and $y^k$. The iterations are then
\begin{subequations}\label{eq:alg_vec}
	\begin{align}
	\xi^{k+1} &= A(L) \xi^k + B(L) u^k \\
	y^k     &= C(L) \xi^k + D(L) u^k \\
	u^k     &= \df(y^k)
	\end{align}
\end{subequations}
where the $i\textsuperscript{th}$ component of $\df(y^k)$ is $\df_i(y_i^k)$, and the system matrices are:
\begin{align*}
A(L) &\defeq I_n\otimes A_0 + L\otimes A_1, & \!\!
B(L) &\defeq I_n\otimes B_0 + L\otimes B_1, \\
C(L) &\defeq I_n\otimes C_0 + L\otimes C_1, & \!\!
D(L) &\defeq I_n\otimes D_0 + L\otimes D_1.
\end{align*}
By Assumption~\ref{assumption:graph}, we may write $L = V \Lambda V^\tp$ where $V=\bmat{v_1 & \cdots & v_n}\in\R^{n\times n}$ is an orthogonal matrix of eigenvectors and $\Lambda=\diag(\lambda_1,\ldots,\lambda_n)$ is a diagonal matrix of eigenvalues with $0 = \lambda_1 < \lambda_2 \le \dots \le \lambda_n$.  For $\ell=1,\dots,n$ we define the state, input, and output in direction $v_\ell$ as follows:
\begin{equation}\label{eq:coord_transform}
	\bar{\xi}_\ell^k \defeq (v_\ell^\tp\otimes I_2) \xi^k \qquad
	\bar{u}_\ell^k \defeq v_\ell^\tp u^k \qquad
	\bar{y}_\ell^k \defeq v_\ell^\tp y^k
\end{equation}
In these new coordinates, \eqref{eq:alg_vec} is equivalent to the $n$ decoupled single-input single-output systems
\begin{subequations}\label{eq:sep_sys}
	\begin{align}
	\bar{\xi}_\ell^{k+1} &= (A_0+\lambda_\ell A_1) \bar{\xi}_\ell^k + (B_0+\lambda_\ell B_1) \bar{u}_\ell^k \\
	\bar{y}_\ell^k &= (C_0+\lambda_\ell C_1) \bar{\xi}_\ell^k + (D_0+\lambda_\ell D_1) \bar{u}_\ell^k \\
	\bar{u}_\ell^k &= v_\ell^\tp \df\bigl( v_1\, \bar{y}_1^k + \cdots + v_n\, \bar{y}_n^k \bigr)
	\end{align}
\end{subequations}
for all $\ell \in\{1,\ldots,n\}$ where $\bar{\xi}_\ell^0\in\R^2$. We refer to~\eqref{eq:sep_sys} as the \emph{separated system}, and the corresponding transfer functions are given by
\begin{align}\label{eq:tf_sep}
G_{\lambda_\ell}(z) &=
(C_0+\lambda_\ell C_1) \bigl( z I - (A_0+\lambda_\ell A_1)\bigr)^{-1} (B_0+\lambda_\ell B_1)\nonumber \\ 
& \qquad +(D_0+D_1\lambda_\ell).
\end{align}
Note that $\ell$ indexes the $n$ directions in the eigenspace of~$L$ as opposed to the $n$ agents, which are indexed by $i$.

Given the generic transfer function $G_{\lambda_\ell}(z)$, imposing that algorithm~\eqref{eq:alg} has an optimal fixed point leads to properties that the transfer functions must satisfy.  Here, an optimal fixed point is any point $(\xi_i^\star,y_i^\star,u_i^\star)$ such that $y_i^\star=x^\star$ solves~\eqref{eq:distrop-problem}. We summarize these in the following lemma; see Appendix~\ref{sec:tfproof} for a proof.

\begin{lem}\label{lem:tf}
	Suppose algorithm~\eqref{eq:alg} converges to an optimal fixed point.  Then the transfer function~\eqref{eq:tf_sep} satisfies the following properties.
	\begin{enumerate}
		\item $G_{\lambda_1}(z)$ has a pole at $z=1$ and is marginally stable.
		\item $G_{\lambda_\ell}(z)$ has a zero at $z=1$ and is strictly stable for all $\ell = 2,\dots,n$.
	\end{enumerate}
\end{lem}

The requirements on the transfer function in Lemma~\ref{lem:tf} imply that no single-state algorithm of the form~\eqref{eq:alg} can solve the distributed optimization problem. This provides an explanation as to why the Distributed Gradient Descent algorithm must use a diminishing stepsize~\cite{DGD}. We therefore restrict ourselves to algorithms with two states (i.e., $\xi_i^k\in\R^{2d}$) since these are the simplest algorithms that can achieve linear convergence rates.

\begin{cor}\label{cor:onestate}
No algorithm satisfying properties~\ref{P2} and~\ref{P3} where each agent has a local state of dimension~$d$ can solve the distributed optimization problem~\eqref{eq:distrop-problem}.
\end{cor}

\paragraph{Proof of Corollary~\ref{cor:onestate}.}
Such an algorithm can be written in the form~\eqref{eq:alg} with $A_i,B_i,C_i,D_i\in\R$ for $i=1,2$ (since the state is dimension $d$) and $D_0=D_1=0$. The corresponding transfer function then has the form
\begin{align*}
G_{\lambda_{\ell}}(z) = \frac{(C_0+\lambda_{\ell} C_1) (B_0+\lambda_{\ell} B_1)}{z-(A_0+\lambda_{\ell} A_1)}.
\end{align*}
Suppose the algorithm solves the distributed optimization problem~\eqref{eq:distrop-problem}. Then from Lemma~\ref{lem:tf}, $G_{\lambda_0}(z)$ must have a pole at $z=1$, and $G_{\lambda_\ell}(z)$ must have a zero at $z=1$ for all $\ell=2,\dots,n$. The first condition implies that $A_0=1$ with $B_0 C_0 \neq 0$, while the second condition implies that either $B_0=B_1=0$ or $C_0=C_1=0$. These conditions contradict each other, implying that the algorithm does \emph{not} solve the distributed optimization problem. \qedhere


\section{Putting algorithms in canonical form}\label{sec:algorithms}

We now outline a procedure to put any distributed optimization algorithm in our class of interest in canonical form. Note that the canonical form has the form~\eqref{eq:alg} with
\begin{align}\label{eq:realization}
\begin{aligned}
\stsp{A_0}{B_0}{C_0}{D_0} &= \left[\begin{array}{cc|c} 1 & \zeta_0 & -\alpha \\ 0 & 1 & 0 \\\hlinet 1 & 0 & 0 \end{array}\right]\quad\text{and} \\
\stsp{A_1}{B_1}{C_1}{D_1} &= \left[\begin{array}{cc|c} -\zeta_1 & \zeta_2 & 0 \\ -1 & 0 & 0 \\\hlinet -\zeta_3 & 0 & 0 \end{array}\right]
\end{aligned}
\end{align}
and the corresponding transfer function~\eqref{eq:tf_sep} is
\begin{align}\label{eq:tf_canonical}
G^\text{CF}_{\lambda_{\ell}}(z) = -\frac{\alpha\,(1-\zeta_3\lambda_\ell) (z-1)}{(z-1)(z-1+\zeta_1 \lambda_{\ell}) + \lambda_{\ell}\,(\zeta_0 + \zeta_2 \lambda_{\ell})}.
\end{align}
We can then apply the following procedure to put any distributed algorithm satisfying properties \ref{P1}--\ref{P3} into our canonical form. We summarize the results for several known algorithms in Table~\ref{table::comparison}.

\begin{enumerate}[itemsep=-0.8mm,topsep=0.8mm]
\item Write the algorithm in the form~\eqref{eq:alg} to obtain the system matrices $(A_0,B_0,C_0,D_0)$ and $(A_1,B_1,C_1,D_1)$.
\item Compute the transfer function $G_{\lambda_{\ell}}(z)$ using~\eqref{eq:tf_sep}.
\item Obtain parameters $\alpha,\zeta_0,\zeta_1,\zeta_2,\zeta_3$ by comparing coefficients of the transfer function $G_{\lambda_{\ell}}(z)$ with that of the canonical form $G^\text{CF}_{\lambda_{\ell}}(z)$ in~\eqref{eq:tf_canonical}.
\end{enumerate}
\begin{table}[htb]
\caption{Parameters of distributed optimization algorithms in canonical form using $W = I-L$.}\label{table::comparison}
\centering
\renewcommand{\arraystretch}{1.1}%
\begin{tabular}{p{3.7cm}cccc}
	& $\zeta_0$ & $\zeta_1$ & $\zeta_2$ & $\zeta_3$ \\ \hline
	EXTRA~\cite{EXTRA}
	& $\frac{1}{2}$          & $1$       & $0$                    & $0$ \\
	NIDS~\cite{NIDS}
	& $\frac{1}{2}$          & $1$       & $0$                    & $\frac{1}{2}$ \\
	Exact Diffusion~\cite{ExactDiffusion}
	& $\frac{1}{2}$          & $1$       & $0$                    & $\frac{1}{2}$ \\
	DIGing~\cite{DIGing,QuLi}
	& $0$              & $2$      & $1$                & $0$ \\
	AsynDGM\tablefootnote{The similar algorithm AugDGM~\cite{AugDGM} does not fit in our framework because it requires two sequential rounds of communication per iteration.}~\cite{AsynDGM}
	& $0$              & $2$      & $1$                & $1$ \\
	Jakoveti\'{c}~\cite{Unification} ($\mathcal{B}=\beta I$)
	& $\alpha\beta$ & $2$      & $1$                & $0$  \\
	Jakoveti\'{c}~\cite{Unification} ($\mathcal{B}=\beta W$)
	& $\alpha\beta$ & $2$      & $1-\alpha\beta$ & $0$ \\
\end{tabular}\vspace{-6mm}
\end{table}

\noindent Since the parameters of our canonical form are unique, we can compare various algorithms by putting them in this form. For example, from Table~\ref{table::comparison} we observe that NIDS and Exact Diffusion have the same parameters and are therefore equivalent. Similarly, for the algorithm of Jakoveti\'{c}, choosing  $\mathcal{B}=0$ recovers DIGing, and $\mathcal{B}=\frac{1}{\alpha} W$ with $W\mapsto\frac{1}{2}(I+W)$ recovers EXTRA, as noted in~\cite{Unification}.

\begin{rem}
To provide an additional degree of freedom, we can relate the gossip matrix to the Laplacian matrix by $W = I-\mu L$ where $\mu\in\R$ is an additional scalar parameter. This is equivalent to choosing the gossip matrix as a convex combination of $I-L$ and the identity, i.e., $W = (1-\mu) I + \mu (I-L)$.
\end{rem}

\begin{rem}[Uniqueness of realization]\label{rem:non_unique_realization}
The canonical form~\eqref{eq:realization} can be reliably obtained from an algorithm by computing its transfer function~\eqref{eq:tf_canonical} and then extracting coefficients. In contrast, working with realizations directly can be problematic because even if coordinates are chosen such that $(A_0,B_0,C_0,D_0)$ has the desired form, this does not ensure uniqueness of $(A_1,B_1,C_1,D_1)$. For example, one can easily check that the following realization has the same transfer function as~\eqref{eq:realization}:
\begin{align*}
\stsp{A_0}{B_0}{C_0}{D_0} &= \left[\begin{array}{cc|c} 1 & \zeta_0 & -\alpha \\ 0 & 1 & 0 \\\hlinet 1 & 0 & 0 \end{array}\right] \quad\text{and} \\
\stsp{A_1}{B_1}{C_1}{D_1} &= \left[\begin{array}{cc|c} -\zeta_1 & \zeta_2 & \alpha\zeta_3 \\ -1 & 0 & 0 \\\hlinet 0 & 0 & 0 \end{array}\right].
\end{align*}
The equivalence class of realizations for a doubly-indexed transfer function (i.e. with $z$ and $\lambda$) is more involved than the singly-indexed case and a broader discussion on this topic is outside the scope of the present work; seminal references include~\cite{fornasini,roesser}.
\end{rem}

\subsection{Example: NIDS}
To illustrate our approach, we return to the NIDS algorithm from Section~\ref{sec:intro}. As suggested by the authors of~\cite{NIDS}, we choose $\widetilde{W}\defeq (I+W)/2$ where the gossip matrix~$W$ is related to the Laplacian matrix~$L$ by $W = I-L$.

First, we write NIDS in the form~\eqref{eq:alg} as follows:
\begin{equation*}
\begin{aligned}
	\bmat{x^{k+2} \\ x^{k+1} \\ \df(y^k)} &= \bmat{2I-L & \!\frac{1}{2} L-I & \!\alpha(I-\frac{1}{2} L) \\ I & 0 & 0 \\ 0 & 0 & 0} \!\! \bmat{x^{k+1} \\ x^k \\ \df(y^{k-1})} \\
	&\qquad + \bmat{\alpha(\frac{1}{2} L - I) \\ 0 \\ I} \df(y^k) \\
	y^k &= \bmat{I & 0 & 0} \bmat{x^{k+1} \\ x^k \\ \df(y^{k-1})}
\end{aligned}
\end{equation*}
This is equivalent to~\eqref{eq:alg} with state-space matrices
\begin{align*}
&\stsp{A_0}{B_0}{C_0}{D_0} = \left[\begin{array}{ccc|c} 2 & -1 & \alpha & -\alpha \\ 1 & 0 & 0 & 0 \\ 0 & 0 & 0 & 1 \\\hlinet 1 & 0 & 0 & 0 \end{array}\right] \qquad\text{and} \\
&\stsp{A_1}{B_1}{C_1}{D_1} = \left[\begin{array}{ccc|c} -1 & \frac{1}{2} & -\frac{\alpha}{2} & \frac{\alpha}{2} \\ 0 & 0 & 0 & 0 \\ 0 & 0 & 0 & 0 \\\hlinet 0 & 0 & 0 & 0 \end{array}\right],
\end{align*}
which have associated transfer function~\eqref{eq:tf_sep} given by
\begin{align*}
G_{\lambda_{\ell}}(z) = -\frac{\alpha\,(1-\frac{1}{2}\lambda_{\ell})(z-1)}{(z-1) (z-1+\lambda_{\ell}) + \frac{1}{2}\lambda_{\ell}}.
\end{align*}
Comparing coefficients with that of $G_{\lambda_{\ell}}^\text{CF}(z)$ in~\eqref{eq:tf_canonical}, we find that the parameters of the canonical form are $(\zeta_0,\zeta_1,\zeta_2,\zeta_3) = (\frac{1}{2},1,0,\frac{1}{2})$.


\section{Conclusion}\label{sec:conclusion}

In this paper, we derived a canonical form for a large class of first-order distributed algorithms that fulfills existence and uniqueness properties, and features a minimal parameterization. We provided sufficient conditions for an algorithm in canonical form to have a fixed point corresponding to the solution of the distributed optimization problem, as well as necessary conditions for convergence to such a fixed point.
Combined with analysis tools such as those described in~\cite{sundararajan_allerton,lessard16,iqcadmm_ICML} and references therein, this work provides a first step toward a principled and automated methodology for the analysis and design of distributed optimization algorithms.

\bibliographystyle{abbrv}
{\small \bibliography{canform}}

\section{Appendix}\label{appendix}

\subsection{Proof of Theorem~\ref{thm:fixed_point}}\label{sec:fixedpointproof}

\paragraph{All fixed points are optimal.}
Consider a fixed point of the algorithm, given by $(v_{1i}^\star,v_{2i}^\star,y_i^\star,u_i^\star,x_i^\star,w_i^\star)$ for ${i\in\{1,\dots,n\}}$. Applying (C.6) and (C.1) to the fixed point, we have $0 = v_{1i}^\star = \sum_{j=1}^n L_{ij}\,x_j^\star$. This implies $y_i^\star = x_i^\star$ from (C.3), and $x_i^\star = x_j^\star$ for all $i,j\in\{1,\ldots,n\}$ since the graph is connected by Assumption~\ref{assumption:graph}. This shows that the agents reach consensus at the fixed point; we have left to show that the consensus variable is the solution to~\eqref{eq:distrop-problem}. To do so, we sum (C.5) to obtain
\begin{align}
\!\alpha \sum_{i=1}^n u_i^\star = \zeta_0 \sum_{i=1}^n w_i^\star - \zeta_1\! \sum_{i,j=1}^n L_{ij} x_j^\star + \zeta_2\! \sum_{i,j=1}^{n} L_{ij} w_j^\star. \label{eq:fixedpoint_sum}
\end{align}
The last two terms on the right side are both zero since the Laplacian matrix is symmetric by Assumption~\ref{assumption:graph}, and therefore satisfies $L^\tp\1 = 0$. Furthermore, the first term on the right side is zero due to T3. To see this, we sum (C.6) to obtain $\sum_{i=1}^n w_i^{k+1} = \sum_{i=1}^n w_i^k$, so either $\zeta_0=0$ or $\sum_{i=1}^n w_i^\star=0$. Finally, $\alpha\neq 0$ from T1, so we must have $0 = \sum_{i=1}^n u_i^\star = \sum_{i=1}^n \df_i(y_i^\star)$. The fixed point satisfies~\eqref{eq:distrop-problem} and is therefore optimal.

\paragraph{An optimal fixed point exists.}
Define $x^\star$ as in~\eqref{eq:distrop-problem}. We show how to construct an optimal fixed point of the canonical form. Define ${v_{1i}^\star \defeq 0}$, ${y_i^\star \defeq x^\star}$, ${u_i^\star \defeq \df_i(x^\star)}$, and ${x_i^\star \defeq x^\star}$ for all ${i\in\{1,\ldots,n\}}$. We have $\sum_{i=1}^n u_i^\star = 0$ from the definition of $x^\star$, so from~\ref{T2}, the linear system
\begin{align}
\alpha u_i^\star = \zeta_0 w_i^\star + \zeta_2 \sum_{j=1}^{n} L_{ij}\,w_j^\star \label{eq:lineq}
\end{align}
has a solution $w_i^\star$. Finally, we define $v_{2i}^\star \defeq \sum_{j=1}^n L_{ij}\,w_i^\star$. Then the point $(v_{1i}^\star,v_{2i}^\star,y_i^\star,u_i^\star,x_i^\star,w_i^\star)$ for $i\in\{1,\dots,n\}$ is an optimal fixed point. \qedhere

\subsection{Proof of Theorem~\ref{thm:realization}}\label{sec:realizationproof}

Consider an algorithm that satisfies properties \ref{P1}--\ref{P3} and converges to an optimal fixed point. Such an algorithm can be represented by the iterations~\eqref{eq:alg} with the following: i) $A_i\in\R^{2\times 2}$ for $i=1,2$ (since the local state vector is in $\R^{2d}$, ii) $D_0=D_1=0$, and iii) either $B_1=0$ or $C_1=0$ (otherwise the algorithm would require two sequential rounds of communication per iteration). We can therefore parameterize the state-space matrices as
\begin{align*}
\addtolength{\arraycolsep}{-2pt}\stsp{A_0}{B_0}{C_0}{D_0} &\!=\! \left[\addtolength{\arraycolsep}{-2pt}\begin{array}{cc|c} a_1 & a_2 & b_1 \\ a_3 & a_4 & b_2  \\\hlinet c_1 & c_2 & 0 \end{array}\right]\!,\!\!\! &
\addtolength{\arraycolsep}{-2pt}\stsp{A_1}{B_1}{C_1}{D_1} &\!=\! \left[\addtolength{\arraycolsep}{-2pt}\begin{array}{cc|c} a_5 & a_6 & b_3 \\ a_7 & a_8 & b_4  \\\hlinet c_3 & c_4 & 0 \end{array}\right]\!\!.
\end{align*}
The corresponding transfer function has the form
\begin{align*}
G_\lambda(z) = \frac{(\eta_1 + \eta_2\lambda + \eta_3\lambda^2) z + (\eta_4 + \eta_5\lambda + \eta_6\lambda^2 + \eta_7\lambda^3)}{z^2 + (\eta_8 + \eta_9\lambda) z + (\eta_{10} + \eta_{11} \lambda + \eta_{12}\lambda^2)}
\end{align*}
where the parameters $\{\eta_i\}$ are defined as follows:
\begin{align*}
\eta_1    &\defeq b_1 c_1 + b_2 c_2 \\
\eta_2    &\defeq b_1 c_3 + b_3 c_1 + b_2 c_4 + b_4 c_2 \\
\eta_3    &\defeq b_3 c_3 + b_4 c_4 \\
\eta_4    &\defeq -a_1 b_2 c_2 + a_2 b_2 c_1 + a_3 b_1 c_2 - a_4 b_1 c_1 \\
\eta_5    &\defeq a_2 b_2 c_3 - a_1 b_4 c_2 - a_1 b_2 c_4 + a_2 b_4 c_1 \\
          &\qquad + a_3 b_1 c_4 + a_3 b_3 c_2 - a_4 b_1 c_3 - a_4 b_3 c_1 \\
          &\qquad - a_5 b_2 c_2 + a_6 b_2 c_1 + a_7 b_1 c_2 - a_8 b_1 c_1 \\
\eta_6    &\defeq  a_2 b_4 c_3 - a_1 b_4 c_4 + a_3 b_3 c_4 - a_4 b_3 c_3 \\
          &\qquad - a_5 b_2 c_4 - a_5 b_4 c_2 + a_6 b_2 c_3 + a_6 b_4 c_1 \\
          &\qquad + a_7 b_1 c_4 + a_7 b_3 c_2 - a_8 b_1 c_3 - a_8 b_3 c_1 \\
\eta_7    &\defeq a_6 b_4 c_3 - a_5 b_4 c_4 + a_7 b_3 c_4 - a_8 b_3 c_3 \\
\eta_8    &\defeq -(a_1+a_4) \\
\eta_9    &\defeq -(a_5+a_8) \\
\eta_{10} &\defeq a_1 a_4 - a_2 a_3 \\
\eta_{11} &\defeq a_1 a_8 - a_2 a_7 - a_3 a_6 + a_4 a_5 \\
\eta_{12} &\defeq a_5 a_8 - a_6 a_7
\end{align*}
Since either ${B_1=0}$ or ${C_1=0}$, we have $\eta_3=\eta_7=0$. Since the algorithm converges to an optimal fixed point, we have from Lemma~\ref{lem:tf} that $G_\lambda(z)$ must have a zero at $z=1$ for all nonzero eigenvalues $\lambda$ of $L$, so the term $(z-1)$ must factor from the numerator. Therefore, the parameters satisfy $\eta_1 + \eta_4 = 0$, $\eta_2 + \eta_5 = 0$, and $\eta_6 = 0$. Also, the transfer function must have a pole at $z=1$ when $\lambda=0$, which requires the denominator to be $(z-1)^2$ when $\lambda=0$. This implies $\eta_8=-2$ and $\eta_{10}=1$. The transfer function then has the form
\begin{align}
G_\lambda(z) = \frac{(\eta_1 + \eta_2\lambda)(z-1)}{(z-1)^2 + \lambda\,(\eta_{11} + \eta_9 z + \eta_{12}\lambda)},\label{eq:tf_eta}
\end{align}
which is equivalent to~\eqref{eq:tf_canonical} with $(\eta_1,\eta_2,\eta_9,\eta_{11},\eta_{12}) \mapsto (-\alpha,\alpha\zeta_3,\zeta_1,\zeta_0-\zeta_1,\zeta_2)$. Note that $\eta_1$ cannot be zero, since this would violate the necessary condition that the transfer function has a pole at $z=1$ when $\lambda=0$. We can then invert the mapping to obtain the parameters $(\alpha,\zeta_0,\zeta_1,\zeta_2,\zeta_3) = \bigl(-\eta_1,\eta_9+\eta_{11},\eta_9,\eta_{12},-\frac{\eta_2}{\eta_1}\bigr)$. Therefore, the algorithm can be put into canonical form with a suitable choice of parameters $\alpha,\zeta_0,\zeta_1,\zeta_2,\zeta_3$, and these parameters are unique since the mapping to the transfer function coefficients is one-to-one. Furthermore, there are five degrees of freedom in the coefficients of the transfer function~\eqref{eq:tf_eta}, so any fully expressive canonical form must have at least five parameters. Our canonical form has precisely five parameters and is therefore a minimal parameterization. It remains to show that the technical conditions hold.

\paragraph{T1.} Suppose $\alpha=0$. Then the transfer function of the canonical form~\eqref{eq:tf_canonical} is identically zero, and therefore does not satisfy the conditions of Lemma~\ref{lem:tf}. But this contradicts that the algorithm converges to an optimal fixed point, so $\alpha\neq 0$. Therefore, \ref{T1} holds.

\paragraph{T2.} For a given initial state, let $(v_{1i}^\star,v_{2i}^\star,y_i^\star,u_i^\star,x_i^\star,w_i^\star)$ for $i\in\{1,\dots,n\}$ denote the optimal fixed point to which the algorithm converges. Then this point must satisfy the linear equation~\eqref{eq:lineq}. For this to hold for a general objective function, the system must have a solution for any $u_i^\star$ such that $\sum_{i=1}^n u_i^\star = 0$, which implies \ref{T2} holds.

\paragraph{T3.} Since the algorithm converges to an optimal fixed point, the right side of~\eqref{eq:fixedpoint_sum} must be zero, or equivalently, $\zeta_0 \sum_{i=1}^n w_i^\star = 0$. The Laplacian matrix is symmetric by Assumption~\ref{assumption:graph}, so $\sum_{i=1}^n w_i^{k+1} = \sum_{i=1}^n w_i^k$ for all $k\geq 0$. Therefore, we must have either $\zeta_0=0$ or $\sum_{i=1}^n w_i^0=0$, so \ref{T3} holds. \qedhere

\subsection{Proof of Lemma~\ref{lem:tf}}\label{sec:tfproof}

Since algorithm~\eqref{eq:alg} converges to an optimal fixed point, the separated systems~\eqref{eq:sep_sys} also converge to fixed points. Since $\lambda_1=0$ and $v_1 = \tfrac{1}{\sqrt{n}}\1$, and by~\eqref{eq:coord_transform} and~\eqref{eq:distrop-problem},
\vspace{-2mm}
\begin{subequations}\label{eq:opt_cond_sep}
	\begin{equation}\label{eq:opt_cond_sepa}
	\bar u_1^\star = v_1^\tp u^\star
	= \tfrac{1}{\sqrt n}\sum_{i=1}^n \df_i(y_i^\star)
	= \sqrt{n}\, \df(x^\star) = 0.
	\end{equation}
	For $\ell=2,\dots,n$, each $v_\ell$ is orthogonal to $v_1$ since $L$ is symmetric. Again using~\eqref{eq:coord_transform} and~\eqref{eq:distrop-problem},
	\begin{equation}\label{eq:opt_cond_sepb}
	\bar y_\ell^\star = v_\ell^\tp y^\star = v_\ell^\tp \1 x^\star = 0
	\qquad\text{for }\ell \in\{2,\dots,n\}.
	\end{equation}
\end{subequations}
We first prove that $G_{\lambda_1}(z)$ has a pole at $z=1$ and is marginally stable.  Suppose by contradiction that all poles of $G_{\lambda_1}(z)$ are strictly stable.  Then for every input sequence, $\bar u_1^k \rightarrow \bar u_1^\star=0$ by~\eqref{eq:opt_cond_sepa}, and the corresponding output sequence $\bar y_1^k$ must tend to $\bar y_1^\star=0$.  However, this is a contradiction because $\bar y_1^\star$ is nonzero in general.  Alternatively suppose that $G_{\lambda_1}(z)$ has an unstable pole.  Then there exists an input sequence such that $\bar y_1^k \rightarrow \infty$ as $\bar u_1^k\rightarrow u^\star$.  This contradicts that $\bar y_1^\star$ is finite.  Hence, $G_{\lambda_1}(z)$ is marginally stable with a pole at $z=1$.  

To prove that $G_{\lambda_\ell}(z)$ is strictly stable for $\ell=2,\dots,n$, assume the contrary.  Then there exists a pole for some $z$ such that $|z|\geq 1$.  Consequently, for a converging input sequence, $\bar y_\ell^k \rightarrow \bar y_\ell^\star \ne 0$, which contradicts~\eqref{eq:opt_cond_sepb}.  Finally, we show that $G_{\lambda_\ell}(z)$ has a zero at $z=1$.  Let $\bar u_\ell^k\rightarrow \bar u_\ell^\star$, which in general is a nonzero constant.  Then the steady-state gain to $y_\ell$ must be zero to ensure~\eqref{eq:opt_cond_sepb} holds. \qedhere

\end{document}